\documentclass[11pt]{amsart}

\usepackage[english]{babel}
\usepackage[utf8]{inputenc}

\usepackage{amsmath,amsthm, amssymb, bbm}
\usepackage{graphicx}
\usepackage{fullpage}

\newcommand{\R}{{{\mathbb {R}}}}

\let\oldmarginpar\marginpar
\renewcommand\marginpar[1]{\-\oldmarginpar[\raggedleft\scriptsize #1]%
{\raggedright\scriptsize #1}}

\begin{document}

\title{A note on Ising random currents, Ising-FK, loop-soups and the Gaussian free field}
\author{Titus Lupu \and Wendelin Werner}
\address {
Institute for Theoretical Studies,
ETH Z\"urich,
Clausiusstr. 47,
8092 Z\"urich,
Switzerland}
\email
{titus.lupu@its-eth.ethz.ch}

\address{
Department of Mathematics,
ETH Z\"urich,
 R\"amistr. 101,
8092 Z\"urich, Switzerland}
\email
{wendelin.werner@math.ethz.ch}

\begin {abstract}
We make a few elementary observations that relate directly the items mentioned in the title. In particular, we note that when one superimposes the random current model related to the Ising model with
an independent Bernoulli percolation model with well-chosen weights, one obtains exactly the FK-percolation (or random cluster model) associated with the Ising model. We also point out that this relation can be interpreted via 
loop-soups, combining the description of the sign of a Gaussian free field on a discrete graph knowing its square (and the relation of this question with the FK-Ising model) with the loop-soup interpretation of the random current model.
\end {abstract}

\maketitle

\section {A simple direct Ising-random-current/Ising-FK coupling}

Let us first briefly review the definitions of the basic models (Ising, random current and FK-Ising) that we will discuss. 
Troughout this section, we will consider a finite connected graph $G$ consisting of a set of vertices $X$ and a set of non-oriented edges $E$. 
We will also use a function $\beta = ( \beta_e)$ from the set of edges into $\R_+$.

\medbreak \noindent
{\sc The Ising model.}
The Ising model on $G$ with edge-weights $(\beta_e)$ is the probability measure on $\Sigma:= \{ -1, +1 \}^X$ defined by
$$P_\beta ( ( \sigma_x) ) = Z_\beta^{-1} \prod_{e} 
\exp ( \beta_e I_\sigma (e)) , $$
 where $I_\sigma (e)$ denotes the product $\sigma_x\sigma_y$ where $x$ and $y$ are the two extremities of the edge $e$, and 
 $Z_\beta $ is the renormalization constant (sometimes called the partition function) chosen 
so that this is a probability measure.

\medbreak \noindent
{\sc The random current model.} The random current model is closely related to the Ising model, and has been instrumental to prove some of its important proprties (see \cite {Ai}, \cite {ADS} and the references therein). 
Here, one assigns to each edge $e$ of the graph a random non-negative integer $N_e$. In fact, our set ${\mathcal N}$ of admissible configurations imposes the 
additional constraint that for each site $x$, the sum of the $N_e$'s for all adjacent vertices to $e$ is even (when $e$ is a vertex from $x$ to $x$, then $N_e$ is counted twice).
The random current model is the probability measure on ${\mathcal N}$ defined by
$$\hat P_\beta ( (n_e) ) = \hat Z_\beta^{-1} \prod_e \frac {(\beta_e)^{n_e}}{n_e!}.$$ 
 By expanding each $\exp (\beta_e I_e (\sigma))$ in the definition of the Ising model into $\sum_{n_e} (\sigma_x \sigma_y \beta_e)^{n_e} / (n_e)!$, and resumming over all $\sigma$'s (noting that all terms with odd degree in $\sigma_x$ do sum up to $0$ by symmetry), one easily sees that the partition function $\hat Z_\beta$ for the random current model is indeed the same 
as the partition function $Z_\beta$ of the Ising model.

\medbreak \noindent
{\sc The FK-Ising model.} The FK-model (named after Kasteleyn and Fortuin) associated to the Ising model (it is also called the random cluster model but we stick here to FK-Ising percolation terminology in order not to confuse random currents with random clusters) is a probability measure on $\{ 0, 1\}^E$. The edge $e$ is said to be open for the configuration $w=(w_e)$ if $w_e=1$, and when 
$w_e = 0$ it is closed. To each configuration $w$, one look at  the graph $G_w$ consisting of the sites $X$ and the 
set of edges that are open for $w$, and denote by $k(w)$ its number of connected components. The FK-Ising model is defined by
$$ \tilde P_\beta ( (w_e)) = \tilde Z_{\beta}^{-1} \times 2^{k(w_e)} \times 
 \prod_e (1- \exp(-2\beta_e))^{w_e}(\exp ( -2\beta_e) )^{1-w_e} $$
 (this is sometimes described as the FK model for $q=2$ and edge-probabilities $1- e^{-2 \beta_e}$). 
The FK-model is  useful to study the Ising model as the connectivity properties of the FK model correspond to the correlation functions of the Ising model. Indeed, when one samples the FK model and then
chooses in an i.i.d. way a sign for each of the clusters of $G_w$, one obtains a function (assigning a sign to each site) that follows exactly $P_\beta$ (and this leads to a simple expression for the ratio between $Z_\beta$ and $\tilde Z_\beta$), see for instance \cite {GG} for basics about FK-percolation.
One can easily check  that this property (``coloring the FK-clusters at random gives the Ising model'') in fact characterizes the law of the FK-clusters (ie. the information that says which sites are in the same cluster, but not necessarily the information about the state of all edges).

\medbreak \noindent
{\sc Bernoulli percolation.} Bernoulli bond percolation with probabilities $(p_e)$ is the product probability measure on $\{ 0,1 \}^E$, where each edge $e$ is open with respective probability $p_e$. In the sequel, we will use the probabilities $p_e:=  1- \exp ( - \beta_e)$.

\bigbreak
We are now ready to state and prove following  coupling lemma:

\medbreak

\noindent
{{\bf  The ``Current+Bernoulli=FK'' coupling lemma.}}
{\sl Let us consider a random current configuration $(N_e)$ with parameters $(\beta_e)$, and an independent Bernoulli percolation configuration $(\xi_e)$ with probabilities $(p_e= 1- \exp ( - \beta_e))$.
We then define $V_e := 1 - 1_{N_e = \xi_e = 0} \in \{ 0, 1\}$ for each $e$, so that $V_e$ is equal to $1$ if and only if at least one of $N_e$ or $\xi_e$ is non-zero. 
Then, the law of $(V_e)$ is exactly the FK-Ising measure on $G$ with weights $(\beta_e)$.}

\medbreak

Note that this provides a direct coupling between the Ising model and the random current model at the level of probability measures, that somehow enlightens the identity between the  partition functions. 

\medbreak

\noindent
{\bf Proof.}
As one can expect for such a simple statement, the proof is fairly  simple as well: 
Let us first define $U_e$ to be equal to $0$, $1$ or $2$, depending on whether $N_e$ is $0$, odd, or positive even. Note that this is enough information in order to construct $V$, and that the law of $(U_e)$
is the probability measure on ${\mathcal N} \cap \{0,1,2\}^E$ defined by
$$ P ( (u_e)) = Z_\beta^{-1} \times \prod_e ( f_e (u_e)),$$
where $f_e (0) = 1$, $f_e (1) = \sinh \beta_e$ and $f_e (2) = (\cosh \beta_e ) -1$.

Let us define $V_e$ as in the statement of the lemma. 
We define $\alpha_e$ to be $1$ if $N_e$ is even, and $-1$ if $N_e$ is odd. Then, we define $\tilde V_e = \alpha_e V_e$. Note that $\tilde V_e = -1$ if and only if $U_e=1$ (the corresponding weight contribution in the probability is therefore $\sinh (\beta_e)$), and that $\tilde V_e = 1$ if either $u_e=2$, or if $u_e=0=1-\xi_e$ (the total weight contribution in the probability is then $(\cosh \beta_e ) - 1 + (1 - \exp (- \beta_e) ) = \sinh (\beta_e)$ as well -- the fact that these two quantities are equal is the key point in the proof). 
It therefore follows, that the probability that $(\tilde V_e) = (\tilde v_e)$ for an admissible $(\tilde v_e)$  (meaning that each site must have an even number of incoming edges with negative $\tilde v_e$), is equal to
$$ Z_\beta^{-1} \times  \prod_e (( \sinh \beta_e )^{|\tilde v_e|} \times ( \exp (-\beta_e))^{1-| \tilde v_e|}).$$ Hence, the probability that 
 $(V_e)=(v_e)$ is equal to 
$$Z_\beta^{-1} \times K_v \times  \prod_e (( \sinh \beta_e )^{v_e} \times ( \exp (-\beta_e))^{1-v_e})$$ 
where $K_v$ denotes the number of admissible choices for $\tilde v$ that are compatible with $v$ (meaning that $|\tilde v_e | = v_e$). In other words, $K_v$ is the number of ways to assign a sign to each open edge $e$ for the configuration $v$, in such a way that each site has an even number of negative incoming signs. 
But this quantity is
easily shown (see below) to be equal to $2^{o(v) + k(v) - |X|}$, where $o(v)$ is the number of open edges for $v$ and $|X|$ the number of sites in the graph, so that the law of $(V_e)$ is 
$$ P ((v_e)) =  
 ( Z_{\beta}^{-1} \times 2^{-|X|})  \times 2^{k(w)} \times  \prod_e ((2 \sinh (\beta_e))^{v_e}\times (\exp (-\beta_e ))^{1-v_e})  ,$$ which is indeed the same as the FK-Ising measure. 
 
In order to see  that $K_v =  2^{o(v) + k(v) - |X|}$, one can for instance proceed by induction, adding edges one-by-one to a forest-like graph (if $G_v$ is a forest, all of its edges have to be of positive sign) and to see that for each new edge that one adds to $v$ without joining two connected components, one gets an additional multiplicative factor $2$ (this is classical; the number $o(v)+k(v) - |X|$ is  the first Betti number, also known as Kirkhoff's cyclomatic number, of the graph $G_v$).
\qed

\medbreak

It is worth noticing that this property of the random current measure trace does in fact characterize the distribution of the configuration $(W_e) := ((1_{N_e \not= 0 }))$ that describes what edges are occupied by the random current. More precisely, the distribution of $(W_e)$ is the only one such that if one considers an independent Bernoulli percolation $(\xi_e)$ with parameters $(p_e)$, and looks at the collection $(\max ( \xi_e, W_e))$, one obtains exactly the FK-Ising model with parameters $(\beta_e)$. Indeed, one can recover by induction over $n \ge 0$, the probability of all configurations with $n$ occupied edges (for instance, the probability that all edges are unoccupied for $W$ is the ratio between the probability that they are all closed for the FK model and the probability that they are all closed for the percolation, and then one can work out the probability of a configuration where just given edge is occupied etc.). 

\section {Relation to loop-soup clusters}

Let us now explain how the previous considerations can be embedded in the setting of the coupling between loop-soup clusters and the Gaussian free field (GFF) as pointed out in \cite {Lupu1}, using the relation between 
random currents and loop-soups described in \cite {W}. This will follow by combining the following observations: 

\begin {itemize}
\item
Consider a discrete GFF $h$ on the graph $G$, where we view the $(\beta_e)$'s as conductances of an electric network. This is the Gaussian random vector $(h_x)$ with intensity proportional to 
$$ \exp ( - \sum_{e \in E} \beta_e |\nabla_e h|^2 / 2) \prod_{x \in X} dh_x $$ 
where $|\nabla_e h| := |h(x) - h(y)|$ where $x$ and $y$ are the two extremities of the edge $e$.   
This GFF is in fact only defined ``up to an additive constant'' (ie. it is not well-defined) because the previous quantity is invariant when one adds the same constant to all $h_x$'s, but we can for instance artificially (and arbitrarily for what will follow, because we will then anyway condition this GFF by the value of its square) add an edge to our connected graph,
joining a site $x_0 \in X$ to a boundary site $o$, and add the condition that $h( o) = 0$ which amounts to multiply the previous expression by $\exp ( - h_{x_0}^2)$
 and ensures that it is integrable on $\R^X$.

One can then easily make sense of the GFF conditioned by the values of its square $(h(x)^2)$ ie. by $h(x)^2 = u (x)$ for a given vector $(u(x))$ in $(0,\infty)^X$: The unknown random quantities are then the signs $\sigma_x$ of $h(x)$, and the conditional distribution of $(\sigma_x)$ is just proportional to the corresponding Gaussian densities at $( \sigma_x u(x))$. One can note that for a given $(u(x))$, 
this density is proportional to the product over all edges $e =(x,y)$ of 
$  \exp (  \beta_e^u  \sigma_x \sigma_y ) $, where the modified weights $ J_e^u $ are defined by
$$  \beta_e^u := \beta_e \times  u(x) \times  u(y).$$ In other words,  the conditional distribution of these signs $(\sigma_x)$ given  $(h(x)^2) = (u(x))$ is exactly an Ising model with weight $(\beta_e^u)$ on the graph $G$. For instance, in the case where one conditions $h(x)^2$ to be equal to $1$ at each site, one gets exactly the Ising model with edge-weight function $(\beta_e)$.
As explained before, one way to sample this is to choose the signs by tossing independent fair coins for each of the clusters of an FK-model with parameters $(\beta_e)$.   

\item 
We now are recall the relation between the square of the GFF and loop-soups on the graph (see \cite {LJ}): The squared GFF is the cumulated occupation time of a continuous-time loop-soup defined on the discrete graph $G$ (where one adds the boundary point $o$ where the process is killed). 
As noted in \cite {Lupu1}, these continuous-time loop soups can be also viewed (see \cite {Lupu1}) as the trace on the sites of a Brownian loop-soup defined on the cable-system associated with the graph (each edge is replaced by a one-dimensional segment on which the Brownian motion can move continuously) -- the time spent by the discrete loops at sites corresponds to the local time spent at this site by the corresponding 
Brownian loop. This provides a natural coupling between the GFF on the discrete graph, the GFF on the cable system, the continuous-time loop-soup on the discrete 
graph, and the continuous loop-soup on the cable system, that we will from now on always implicitely use. A key observation in \cite {Lupu1} is that conditionally on this loop-soup (that defines the square of the GFF), the sign of the GFF is chosen to be constant and independent for each ``cluster'' of the cable system loop-soup.

\item  
But it is also easy to make sense of the distribution of the loop-soups (on the discrete graph and on the cable systems) conditioned by the square of the GFF on the sites. Indeed, as explained in 
\cite{W}, the conditional distributions of the number of jumps of the loop-soup on $G$ along the unoriented edges of $G$, when conditioned by $(h^2(x))=(u(x))$ is exactly the random current model with edge-weights $(\beta_e^u)$. For instance, when one conditions $h^2 (x)$ to be equal to $1$ at each site, one gets exactly the random current model on $G$ with edge-weights $(\beta_e)$. 
In the loop-soup on the cable system, it is easy to see (this type of observations is already present in \cite {Lupu1}) that on top of the excursions made by the loops inbetween different sites (these correspond to the discrete jumps that we just described via the random current distribution), one adds an independent contribution in each edge $e$ (these correspond to the excursions away from 
the two extremities of the edge that do not cross the edge, and the loops that are totally contained in this edge). When occupation time of the loop-soup at both extremities of $e$ is equal to one, the condtional probability that these contributions join them into the same cluster is equal $p_e$. 
Hence, the conditional distribution (given that $h(x)^2 =1$ at all sites) of the clusters created by the loop-soup on the cable-system is exactly given by the clusters defined by superimposing of a random current (given by (1)) and a Bernoulli percolation on the edges. Comparing this with the previous description of the conditional distribution of $\sigma$, we conclude that the FK-clusters are indeed distributed like the clusters of the superposition of the random current with the Bernoulli percolation. 
\end {itemize}

This therefore provides an alternative explanation to the relation between the Ising random current and the FK-Ising + Bernoulli percolation pointed out in the previous section. In fact, it is this interpretation of the random current in terms of loop-soups conditioned by the values of the GFF at sites that did lead us to realize that the relation derived in the first section should hold (and then, once one guesses that this result holds, it is actually easy to find a direct proof).

\medbreak

Let us note that the notions of loop-soup clusters and their relation with GFF have a nice SLE/CLE type properties in the two-dimensional continuous space via the Brownian loop-soup introduced in \cite {LW}, see \cite {ShW,Lupu2,QW} and the references therein.

\medbreak

To conclude, let us note that is quite possible that some of these random current-loop soup-FK features have been observed before (explicitely or in some slightly hidden way) -- the study of the Ising model has proved to be prone to recurrent rediscoveries of such simple combinatorial identities...  As in \cite {CL,W}, the present considerations are reminiscent of some ideas in \cite {BFS,Dy,Sy,LJ}.

In the recent work of Sabot and Tarr\`es \cite{SabotTarres}  on vertex reinforced jump processes and Ray-Knight theorems, one can for instance find some traces of the relation between loop-soups, fields and the Ising model. More precisely, one can interpret their ``magnetized inversed VRJP'' as a reconstruction 
of the loop-soup conditioned on its occupation field, that is to say on $h^{2}$, that also samples a random current given the edge weights $(\beta_e)$: 
In their setting, edge weights evolve over time, and one first discovers the loops that go through a point $x_1$, then the loops that go through $x_2$ without visiting $x_1$ and so on. Loosely speaking, tracing the loops then progressively eats up the available time at each sites, and the evolving edge-weights represent this remaining available time.  

More precisely, let $x_{1},\dots,x_{k}$ be an arbitrary enumeration of vertices of $G$. One defines by induction over $i \le k$ the processes 
$(\beta_e^{(i)}(t))_{t\geq 0}$ and 
$(X^{(i)}(t))_{t\geq 0}$ as follows. They start from 
 $\beta_e^{(1)}(0):=\beta_e$, $\beta_e^{(i)}(0):=\lim_{t \to \infty}\beta_e^{(i-1)}(t)$ for $2\leq i\leq k$, and  $X^{(i)}(0)=x_{i}$. The dynamics of the edge-weights  $\beta_e^{(i)}(t)$ is described by
$$d\beta_e^{(i)}(t)=-1_{\{e~\text{adjacent to}~X^{(i)}(t)\} } \times \beta_e^{(i)}(t) dt$$
and the dynamics of the jump process $X^{(i)}$ is that when it is at $x$ at time $t$, it jumps to a neighbour $y$ via the edge $e$, with rate
$$\beta_e^{(i)}(t)  \times  {\langle \sigma_{x_{i}}\sigma_{y} \rangle^{(i)}_{t}}
/ {\langle\sigma_{x_{i}}\sigma_{x}\rangle^{(i)}_{t}}, $$
where 
$$
\langle \sigma_{x}\sigma_{y}\rangle^{(i)}_{t}:= E_{(\beta_e^{(i)} (t))} ( \sigma_x \sigma_y)
$$
can be interpreted as the time-evolving two-point functions of the Ising model associated to the time-evolving weights.
It turns out that almost surely,  $X^{(i)}(t) = 
x_{i}$ for all large $t$. Hence, for any edge $e$ adjacent to $x_{i}$, $\lim_{t \to \infty} \beta_e^{(i)}(t)=0$. In particular $\lim_{t \to \infty} \beta_e^{(k)}(t)=0$  for all $e$. 
The  family $(N_e)$ where $N_e$ denotes the total number of jumps across the edge $e$
by the $k$ processes $X^{(i)}$ is then distributed like a random current with weights $(\beta_e)$. For details on this model, see \cite {SabotTarres}.

\bigbreak
\noindent
{\bf Acknowledgements.}
TL acknowledges the support of Dr. Max Rössler, the Walter Haefner Foundation and the ETH Zurich Foundation. 
WW acknowledges the support of the SNF grant SNF-155922. The authors are also part of the NCCR Swissmap of the SNF.


\begin{thebibliography}{99}

\bibitem {Ai}
M. Aizenman. Geometric analysis of $\varphi^4$ fields and Ising models,
Comm. Math. Phys. 86, 1–48, 1982.

\bibitem {ADS}
M. Aizenman, H. Duminil-Copin, V. Sidoravicius.
    Random currents and continuity of Ising model’s spontaneous magnetization,
    Comm. Math. Phys. 334, 719-742, 2015.


\bibitem {BFS}
D. Brydges, J. Fr\"ohlich, T. Spencer.
The random walk representation of classical spin
systems and correlation inequalities, 
Commun. Math. Phys. 83, 123-150, 1982.

\bibitem {CL}
F. Camia and M. Lis. 
Non-backtracking loup soups and statistical mechanics on spin networks, 
preprint, 2015.

\bibitem {Dy}
E.B. Dynkin.
Markov processes as a tool in field theory, J. Funct. Anal. 50, 167-187, 1983.

\bibitem {GG}
G.R. Grimmett. 
The random-cluster model, Springer.

\bibitem {QW}
W. Qian and W. Werner.
Decomposition of Brownian loop-soup clusters,
preprint, 2015.

\bibitem{LW}
G.~F. Lawler and W.~Werner.
\newblock {The Brownian loop soup},
\newblock {Probability Theory and Related Fields} 128, 565-588, 2004.

\bibitem {LJ}
Y. Le Jan.
Markov paths, loops and fields.
L.N. in Math, 2026, 2011.

\bibitem {Lupu1}
T. Lupu.
 From loop clusters and random interlacement to the free field, 
{Ann. Probab.}, to appear.
 
\bibitem {Lupu2}
 T. Lupu.
  Convergence of the two-dimensional random walk loop soup clusters to CLE,
  preprint, 2015.
  
\bibitem {SabotTarres}
C. Sabot and P. Tarres.
 Inverting Ray-Knight identity, 
{Probability Theory and Related Fields}, to appear.  
  
\bibitem{ShW}
S.~Sheffield and W.~Werner.
{Conformal Loop Ensembles: The Markovian characterization and the
  loop-soup construction},
 {Ann. Math.}, 176, 1827--1917, 2012.

\bibitem {Sy}
K. Symanzik.
Euclidean quantum field theory, In: Local quantum theory. Jost, R. (ed.), Academic Press 1969


\bibitem {W1}
W. Werner. 
Topics on the Gaussian Free Field, 
Lecture Notes, 2014.

\bibitem {W}
W. Werner. 
On the spatial Markov property of soups of unoriented and oriented loops,
preprint, 2015.

\end{thebibliography}
\end{document}